\documentclass{svjour2}                    
\smartqed  
\usepackage{graphicx}

\newcommand{\ZZ}{\mathbb{Z}}
\newcommand{\NN}{\mathbb{N}}
\newcommand{\RR}{\mathbb{R}}

\newtheorem{thm}{Theorem}
\newtheorem{defi}{Definition}
\newtheorem{lem}{Lemma}

\newtheorem{cor}{Corollary}
\newtheorem{pro}{Proposition}
\newtheorem{rem}{Remark}

\usepackage{amsmath}
\usepackage{amssymb,latexsym,eucal}

\usepackage{latexsym}
\begin{document}

\title{Invariant Measures and Decay of Correlations for a Class of Ergodic Probabilistic Cellular Automata
}


\author{Coletti, Cristian F        \and
        Tisseur, Pierre 
}


\institute{Coletti, Cristian F. and Tisseur, Pierre \at
              Centro de Matem\'atica, Cogni\c{c}\~ao e 
Computa\c{c}\~ao\\ Universidade Federal do ABC\\
 Santo Andr\'e S.P, Brasil \\
              \email{\{cristian.coletti, pierre.tisseur\}@ufabc.edu.br}           
}

\date{Received: date / Accepted: date}

\maketitle

\begin{abstract}
We give new sufficient ergodicity conditions for two-state probabilistic cellular automata (PCA) of any dimension and any radius. The proof of this result is based on an extended version of the duality concept.
Under these assumptions, in the one dimensional case,  we study  some properties of the unique invariant measure 
and show that it is shift-mixing. 
Also, the decay of correlation is studied in detail. In this sense, the extended concept of duality gives 
exponential decay of correlation and allows to compute explicitily all the constants involved.
\end{abstract}

\section{Introduction}
\renewcommand{\theequation}{1.\arabic{equation}}
\setcounter{equation}{0}

\paragraph
\noindent Probabilistic cellular automata (PCA) are discrete time Markov processes which have been intensely studied 
since at least Stavskaja and Pjatetskii-Shapiro \cite{Stavskaja} (1968). This kind of processes have as state space 
a product space $X=A^G$ where $A$ is any finite set and $G$ is any locally finite and connected graph. 
In this work we will focus our attention on $G=\mathbb{Z}^d$ and $A=\{0, ..., n\}$ for some integer $n \geq 1$. 
We may regard a PCA as an interacting particle system where particles update its states simultaneously and 
independently. Recall that a PCA is ergodic if there exists only one invariant measure $\mu$ and starting from any 
initial measure $\mu_0$ the system converges to $\mu$. 

The aim of this paper is to use  duality principles to study the ergodicity of two-states PCA.
More precisely our work  gives new sufficient ergodicity conditions for the expression of the PCA's local transition probabilities
(see Theorem \ref{EXPR}) and 
show that under these conditions the invariant measure is shift-mixing with exponential decay of correlation.
Relations between the PCA and the dual process (see Lemma \ref{couplage} and Lemma \ref{pro-dense} )
 also allow us to give a very simple expression of 
the constant of the decay of correlation as a function of the radius (of the PCA) and the transition probabilities of the PCA
(see Theorem \ref{decay}). Moreover the proof of  
Theorem \ref{EXPR}  shows in detail how  
to compute the value of the invariant measure on cylinders.
Results about the decay of correlation is an answer to a question raised in \cite{LSS2006}.

The existence of a dual process satisfying the duality equation (see Definition 1 and Liggett \cite{Lig})
  gives useful information (problems in uncountable sets can be reformulated as problems in countable sets) on the PCA but  is not always sufficient to prove that a PCA is ergodic.
In \cite{LSS2006}, Lopez, Sanz and Sobottka introduced an extended concept of duality 
(see Definition 2) and gave general results about ergodicity (see Theorem \ref{MTHM}). 
They used this powerful general theory to give results on  multi-state one-dimensional PCA of radius one 
 and extended previous results about the {\it Domany-Kinzel model} (see \cite{DK} 
for an introduction and \cite {Konno1}  for extensions). Previously, in \cite{Konno2} Konno has given ergodicity conditions for multi-state one-dimensional 
PCA  using  self-duality equations.

Even if, in some cases, the existence of null transition probabilities allows to prove ergodicity of a certain class of 
PCA (see \cite{Konno1} and \cite{Konno2}), it had been conjectured that in the one-dimensional case 
positive noise cellular automata are ergodic. However, P. Gacks, in 2000
, introduced a very complex counterexample (see \cite{Gacs} and \cite{Gray}) for noisy deterministic 
 cellular automata. In that case, the noisy one-dimensional cellular automata does not forget the past and starting from different 
 initial 
 distribution, the PCA may converge to different invariant measures. 
His result can be extended to noisy PCA with positive rates. This conjecture was formulated only in the one-dimensional case since in dimension 2 or higher, it is  
 easier to show the existence of at least two invariant measures.  
 For example the two dimensional Ising model \cite{Gray} or the Toom example (see \cite{Toom1}) that 
exhibit {\it eroder}  properties.
From Theorem \ref{EXPR} there exists a  subclass of attractive PCA (class $\mathcal{C}$) where 
the noisy conjecture is verified ($p(I_r)<1$ implies ergodicity).

In \cite{LSS2006}, the authors explore some ergodic conditions for multi-state PCA.
When the number of states is greater than 2, the conditions of ergodicity are rather restrictive 
in order to be able to give general results. More general ergodicity conditions 
 are  interesting (see \cite{LSS2006}, Section 3.2) but seems to be very complex when the radius grows. 
In this paper, we restrict the study to the 2 states case, which allows to show more easily 
general results for PCA of any radius. 

These sufficient ergodicity conditions can be compared to the Shlosman-Dobrushin condition applied to PCA (see \cite{MS91} and \cite{MS93}). 
In some examples  (see section 3.1)
our sufficient conditions induced by the concept of duality  allow to show ergodicity and decay of correlations where 
the Dobrushin conditions can not be applied.
Moreover, for some classes of ergodic PCA Theorem \ref{decay} gives  greater constants for the   decay of 
spatial correlation.

This paper is organized as follows. In section $2$, we present the basic definitions, notations and some preliminary results. In section $3$, we state the main results, Theorem \ref{EXPR} and Theorem \ref{decay}. We prove Theorem \ref{EXPR} in section $4$. We conclude the paper in section $5$ with the proof of the decay of correlation.
\section{Definitions, notations and preliminary results}
\subsection{Probabilistic Cellular Automata} 

We give a brief description of the theory of PCA. 

Let $A$ be a finite set and $G$ a countable set. Let $X=A^G$ be endowed with the product topology. A probabilistic cellular automata is a discrete time Markov process on the state space $X$.

Let $M(X)$ be the set of probability measures on $X$ and $F(X)$ the set of real functions on $X$ which depend only on a finite number of coordinates of $G$.

The evolution of probabilistic cellular automata is given through their transition operator.

\begin{defi}
A local transition operator is a linear operator
$$
P:F(X) \rightarrow F(X)
$$
such that $Pf \geq 0$ for all $f \geq 0$ and $PI=I,$ where $I : X \rightarrow X$ is the identity function.
\end{defi} 
\begin{defi}
A local transition operator is called independent if
$$
P\left(\varphi \phi \right) = P \left(\varphi \right) P \left(\phi \right)
$$
for all $\varphi, \phi \in F(X)$ such that the finite subsets of G on which they depend do not intersect.
\end{defi}
The independent local transition operators can be defined through the values $p_x(\eta ,k)$ for $x \in G, \eta \in X$ and $k \in A$, as
$$
Pf\left( \eta \right) = \int_X f\left( \sigma \right) \displaystyle \prod_{x \in G} p_x(\eta ,d \sigma (x)), \ \ \ \ \forall f \in F(X),
$$
where for all site $x \in G$, for all $\eta, p_x(\eta ,.)$ is a probability measure on $A$. 
The value $p_x(\eta ,k)$, 
called transition probabilities, represent the probability that the sites $x \in G$ takes the value $k$ in the 
next transition if the present configuration 
of the system is $\eta$. 
For more details see Toom et al. \cite{Toom2}, Maes and Shlosman \cite{MS91} and Lopez and Sanz \cite{LS2000}.

Let $d\ge 1$ be an integer number, $R$ a finite subset of $\ZZ^d$ of cardinality $\vert R\vert$ and 
$f$ a map from $A^{\vert R\vert +1}$ to $[0,1]$. In the particular case $G=\mathbb{Z}^d$ it follows from the discussion above that the 
discrete  time Markov process $\eta_{.}=\{\eta_t(z)\in A : t \in \NN, z\in \ZZ^d\}$ 
whose evolution satisfies

\begin{equation}\nonumber
\mathbb{P}\left[\eta_{t+1}(z)=a | \eta_t(z+i)=b_i, \ \forall \ i \in R\right] = f\left(a,(b_i)_{i \in R}
\right),
\end{equation} 
for all $t\in\NN$ and $z\in \ZZ^d$
is a well defined (discrete time) stochastic process which from now on will be called $d$-dimensional PCA. Here, $\mathbb{P}$ stands for the probability measure on $A^{\mathbb{Z}^d}$ induced by the family of local transition probabilities. Also, let $\mathbb{E}$ be the expectation operator with respect to this probability measure.

 
Let $\mu_0$ be the initial distribution of the PCA.
For any $t\ge 0$, we call $\mu_t$ the distribution of the process at time $t$.
The measure $\mu_t$ is defined on cylinder $u=N(\Lambda,\phi)=\{\xi \in A^{\ZZ^d} : \xi (x) = 
\phi (x) \ \forall \ x \in \Lambda\}$ for some fixed $\phi \in A^{\mathbb{Z}^d}$ and $\Lambda \subset \mathbb{Z}^d,
 \vert \Lambda \vert < \infty$ by
$$
\mu_t(u)=\displaystyle \sum_{v \in \mathcal{C}_t} \mu_0 (v) \mathbb{P}_{\eta_0\in v}\{\eta_t\in u\},
$$  
where $\mathcal{C}_t$ is the family of all cylinders of $X$ on the coordinates of $\Lambda$ 
 ( the finite subset of $\ZZ^d$ used to defined  $u$).

In this paper the notation $\vert \Lambda\vert$ will represent the cardinality of $\Lambda$ when 
$\Lambda$ is a finite subset of $\ZZ^d$. If $U=N(\Xi ,\phi)$ is a cylinder set, 
the notation $\vert U\vert$ will represent the cardinality $\vert\Xi \vert$ of the set $\Xi \subset \ZZ^d$.
In the one dimensional case we adopt the following notation: 
 For any sequence of letters $U=(u_0,\ldots , u_n) \in A^{n+1}$, the set 
 $[U]_s=[u_0\ldots u_n]_s:=\{x\in A^\ZZ\vert x(s)=u_0,\ldots ,x(s+n)=u_n\}$ will be called cylinder 
and $\vert U\vert=n+1$.
\subsection{Two-state Probabilistic Cellular Automata}

In order to simplify the notation we will focus our attention on two-state PCA, that  is to say  PCA 
$\eta_.$ on  $\{0,1\}^{\mathbb{Z}^d}$. For any positive integer $r$, let us define 
$$
I_r := \{i=\left(i_1,\ldots, i_d\right) \in \mathbb{Z}^d: -r \leq i_1, \ldots, i_d \leq r\}.
$$
Define a family of transition probabilities $\{p(J) : J \subset I_r\}$ by
$$
p(J) := \mathbb{P}\{\eta_{t+1}(z)=1\vert \eta_{t}(z+j)=1 : j \in J\}.
$$
Note that any PCA with state space $\{0,1\}^{\ZZ^d}$ is completely characterized by a positive integer number 
$r$ called the radius of the PCA and the set of transition probabilities $\{p(J) : J \subset I_r\}$.
\subsection{The invariant probability measure}
\begin{defi}
Let $T$ be a measure-preserving transformation of a probability space $\left(X,\mathcal{F},\mu \right)$, where 
$\mathcal{F}$ is the $\sigma$-algebra generated by the cylinder sets on $X$. We say that the probability measure 
$\mu$ is $T$-mixing if $\forall \ U, V \in \mathcal{F}$ 
$$
\lim_{n\to\infty}\mu (U\cap T^{-n} V)=\mu (U) \mu (V).
$$
\end{defi}
Since the cylinder sets generate the $\sigma$-algebra $\mathcal{F}$, it follows that the measure $\mu$ is 
$T$-mixing when the last relation is satisfied by any pair of cylinder sets $U$ and $V$ (for more details see \cite{W}).

\subsection{Duality}
The concept of duality is a powerful tool in the theory of interacting particle system. It provides relevant information 
about the evolution of the process under consideration from the study of other simpler process, the dual process. 
The reformulated problems may be more tractable than the original problems and some progress may be achieved. Now we give 
the (classical) definition of duality taken from \cite{Lig}.
\begin{definition}
Let $\eta_.$ and $\zeta_.$ be two Markov processes with state spaces $X$ and $Y$ respectively, and let 
$H\left(\eta,\zeta \right)$ be a bounded measurable function on $X \times Y$. The processes $\eta_.$ and $\zeta_.$  
are said to be dual to one another with respect to $H$ if
\begin{equation} \label{classicduality} \nonumber
\mathbb{E}^{\eta}\left[H\left(\eta_t, \zeta \right)\right] = \mathbb{E}^{\zeta}
\left[H\left(\eta, \zeta_t \right)\right]
\end{equation}
for all $\eta \in X$ and $\zeta \in Y$.
\end{definition}
Unfortunately, it is not true that every process has a dual. Recently, Lopez et al \cite{LSS2006} gave a new notion 
of duality which extends the previous one.
More precisely, they gave the following definition.
\begin{definition} 
Given two discrete time Markov processes, $\eta_t$ with state space $X$ and $\zeta_t$ 
with state space $Y$ 
and $H:X\times Y \rightarrow \mathbb{R}$ and $\mathcal{D}:Y \rightarrow [0,\infty)$ bounded measurable functions, 
 the process  $\eta_.$ and $\zeta_.$ are said  dual to one another with respect to $\left(H,\mathcal{D}\right)$ if
\begin{equation} \label{ExtendedDUAL}
\mathbb{E}_{\eta_0=x}\left[H\left(\eta_1,y\right)\right]=\mathcal{D}(y) \mathbb{E}_{\zeta_0=y}
\left[H\left(x,\zeta_1\right)
\right].
\end{equation}
\end{definition}
\subsection{Duality and sufficient conditions for ergodicity} 

In order to state our results in section 3, we need to give the spirit and some elements of the proof of the following Theorem,
which appears in \cite{LSS2006}.

\begin{thm}\label{MTHM} \cite{LSS2006}
Suppose $\eta_t$ is a Markov process with state space $\bf{X}$ and $\xi_t$ is a markov chain with countable state 
space $\bf{Y}$, which are dual to one another with respect to $(H,\mathcal{D})$. 
If $0 \leq \mathcal{D}(y) < 1$ for all $y \in \bf{Y}$, then 
there exist a stochastic process $\tilde{\xi}_t$ with state space $\tilde{Y}=\bf{Y} \cup \{\mathcal{S}\}$ 
with $\mathcal{S}$ an extra state and a bounded measurable function 
$\tilde{H}:X\times\tilde{Y}\rightarrow\mathbb{R}$ such that $\eta_.$ and 
 $\tilde{\xi}_.$ are dual to one another with respect to $H$. Furthermore, denoting by $\Theta$ the set 
 of all absorbing states of $\xi_.$, if
\begin{itemize}
\item [i)] the set of linear combinations of $\{H(.,y):y\in \bf{Y}\}$ is  dense in  $C\left(X\right)$, the set of 
continuous maps from 
$\bf{X}$ to $\RR$;
\item [ii)] $\mathcal{D}(y) < 1$ for any $y \notin \Theta$, and ${\bf D} := \sup_{y \in \bf{Y} : 
\mathcal{D}(y) < 1} \{\mathcal{D}(y)\} < 1$;
\item [iii)] $H(.,\theta) \equiv c(\theta)$ for all $\theta \in \Theta$ with $\mathcal{D}(\theta)=1$;
\end{itemize}
then $\eta_.$ is ergodic and its unique invariant measure is determined for any $y \in \bf{Y}$ by
\begin{equation}\label{muchap}
\tilde{\mu}(y) = \sum_{\theta \in \Theta \\ d(\theta)=1} c(\theta) \mathbb{P}_{\tilde{\xi}_0=y}
\left[\tilde{\xi}_{\tau}=\theta\right],
\end{equation}
where $\tau$ is the hitting time of $\{\theta \in \Theta:\mathcal{D}(\theta)=1\}\cup \{\mathcal{S}\}$ for 
$\tilde{\xi}_t$ 
and $\hat{\mu}=\lim_{t\to\infty}\hat{\mu_t}$ with 
$$
\hat{\mu}_t(y)=\int_{\bf X} H(x,y)d\mu_t (x).
$$ 
\end{thm}


\noindent {\bf Sketch of the proof.}
First recall that $\tau$ is the hitting time of the dual process $\tilde{\xi_.}$ entering an absorbing state 
$\theta\in\tilde{\Theta}$.
If there exists a dual process  $\tilde{\xi}$ and a function $\tilde{H}$ that satisfies the 
following (classical) duality equation 
\begin{equation}\label{eqdual}
\mathbb{E}_{\eta_0=x}\left[\tilde{H}\left(\eta_1,y\right)\right]=\mathbb{E}_{\tilde{\xi_0}=y}
\left[\tilde{H}\left(x,\tilde{\xi_1}\right)\right],
\end{equation}
it is possible to show that $\hat{\mu}_s(y)=\int_{\bf X} H(x,y)d\mu_s (x)=\mathbb{E}_{\tilde{\xi}_0=y}
[\hat{\mu}(\tilde{\xi}_s)]$.
If $\mathbb{P}\{\tau<\infty\}=1$, it follows that 
$$
\begin{array}{ll}
\lim_{s\to\infty}\hat{\mu}_s(y)
&=\lim_{s\to\infty}\sum_{\theta\in\tilde{\Theta}}\mathbb{E}_{\tilde{\xi}_0=y}[\hat{\mu}(\tilde{\xi}_s)
\vert \tilde{\xi}_t=\theta,\tau\le s]\mathbb{P}_{\tilde{\xi}_0=y}\{\tilde{\xi}_t=\theta,\tau\le s\}\\
&+\lim_{s\to\infty}\mathbb{E}_{\tilde{\xi}_0=y}[\hat{\mu}(\tilde{\xi}_s)
\vert\tau > s]\mathbb{P}_{\tilde{\xi}_0=y}\{\tau > s\}\\
&=\sum_{\theta\in\tilde{\Theta}} \hat{\mu}(\theta)\mathbb{P}_{\tilde{\xi}_0=y}\{\tilde{\xi}_{\tau}=\theta\}.
\end{array}
$$

Finally, when the set of linear combinations of the 
set $\{\tilde{H}(.,y)\vert y\subset\ZZ^{d}\}$ is dense in $C(\bf{X})$ 
(the set of continuous functions from $\bf{X}$ to $\RR$) 
the sequence $(\mu_n)_{n\in\NN}$ converges in the weak* topology. 
Also, the limit measure $\mu$ does not depend on the initial measure $\mu_0$.
\bigskip

Hence, we have seen that the key point is to prove that $\mathbb{P}\{\tau<\infty\}=1$.
One way to show this, is to introduce the new type of duality (see Equation \ref{ExtendedDUAL}). 
If there exists a dual process $\xi_.$ with state space $\bf{Y}$  that verifies 
the new concept of duality (see Equation \ref{ExtendedDUAL}) then 
we can define a standard dual process 
$\tilde{\xi_.}$ with state space  $\tilde{Y}=\bf{Y}\cup\{\mathcal{S}\}$ and such that the set of all 
absorbing states is $\tilde{\Theta}=\Theta\cup \{\mathcal{S}\}$ where $\Theta$ denote the set of all the 
absorbing states of $\xi_.$. 
Here $\mathcal{S}$ is an extra absorbing state and the transition probabilities of $\tilde{\xi_.}$ satisfy 
$$
\mathbb{P}_{\tilde{\xi}_0=\tilde{y}_0}\{\tilde{\xi}_1=\tilde{y}_1\}=\left\{
\begin{array}{ll}
\mathcal{D}(\tilde{y}_0)\mathbb{P}_{\xi_0=\tilde{y_0}}\{\xi_1=\tilde{y}_1\} &,\mbox{ if } \tilde{y}_0,\tilde{y}_1\in \bf{Y}\\
1-\mathcal{D}(\tilde{y}_0)&,\mbox{ if } \tilde{y}_0\in {\bf Y}, \tilde{y}_1=\mathcal{S}\\
1&, \mbox{ if } \tilde{y}_0=\tilde{y}_1=\mathcal{S}.
\end{array}
\right.
$$
Taking $\tilde{H}(x,y)=H(x,y)$ when $y\in \bf{Y}$ and $\tilde{H}(x,\mathcal{S})=0$  
  we obtain that the dual process $\tilde{\xi_.}$ satisfies the standard duality equation 
 \ref{eqdual} and that $\hat{\mu}(\mathcal{S})=0$. 
Note that since ${\bf D}=\sup_{y \in {\bf Y} : \mathcal{D}(y) < 1} 
\{\mathcal{D}(y)\} < 1$,  at each iteration the probability to enter the extra absorbing state 
$\mathcal{S}$ is positive and this implies the following result:

\begin{lem}\label{tau}
Under the conditions of Theorem 1, for all integer $i\ge 1$ one has 
$$
\mathbb{P}(\tau > i)\le {\bf D}^i.
$$
\end{lem}

\noindent {\bf Proof.}
By the Markov property we have that 
$$
\mathbb{P}_{\tilde{\xi_0}=\tilde{y_0}}\{\tau > i\}\leq {\bf D}
\times\mathbb{P}_{\tilde{\xi_0}=\tilde{y_0}}\{\tau > i-1\}.
$$
Then, the result follows by using the mathematical induction principle.

$\hfill\Box$
\medskip

Note that Lemma \ref{tau} implies that  $\mathbb{P}\{\tau<\infty\}=1$ which finishes the proof of Theorem \ref{MTHM}.

\hfill$\Box$

Before stating the main results of this paper, we introduce one more piece of notation: let $^\infty 1^\infty$ denote the all one configuration, i.e. $^\infty 1^\infty=\left(1_{\mathbb{Z}^d}(x)\right)_{x \in \mathbb{Z}^d}$. Analougsly, $^\infty 0^\infty$ denote the all zero configuration.
\section{Main Results and Examples}
A PCA of radius $r$ is called attractive if for any $J\subset I_r$ and $j\in  I_r$ we have 
$p(J\cup\{j\})\ge p(J)$. We consider here the following subclass of attractive PCA.
For any $r\in\NN$, let $\wp(I_r)$ be the set of all subsets of $I_r$.
We say that a two-state PCA of radius $r$ belongs to $\mathcal{C}$ if its transition probabilities satisfy 
$ p(J)=\sum_{J'\subset J}\lambda (J') $ 
where $\lambda$ is some map from $\wp (I_r)\to [0,1)$. 
\medskip
$\mbox{ }\\$
The definition of the class $\mathcal{C}$ is constructive. The following Proposition gives 
sufficient conditions for an attractive PCA to belong to $\mathcal{C}$.  
\begin{pro}\label{cond}
A two-state probabilistic cellular automaton $\eta_.$ belongs to $\mathcal{C}$ if its transition probabilities satisfy the following set of inequalities:

\noindent (a) For any $i \in I_r$,
\small
\begin{eqnarray}
&&p(\{i\}) \ge p(\emptyset). \nonumber
\end{eqnarray}
\normalsize
\noindent (b) For any $1 \leq k \leq |I_r|-1$ and for any $j_0,\ldots ,j_k\in I_r$
\small
\begin{eqnarray*}
p(\{j_0,\ldots ,j_k\}) \ge (-1)^{k}p(\emptyset) - \sum_{n=0}^{k-1} (-1)^{k+1-n}\sum_{\{l_0,\ldots , l_n\}\subset\{j_0,\ldots ,j_k\}} p(\{l_0,\ldots l_n\}).
\end{eqnarray*}
\normalsize
\end{pro}
\begin{thm}\label{EXPR}
Let $\eta_.$ be a two-states d-dimensional probabilistic cellular automaton of radius $r$
 that belongs to $\mathcal{C}$.
If $p(I_r)<1$ then $\eta_.$ is an ergodic PCA and  there exists a dual process $\xi$ which 
satisfy equation  \ref{ExtendedDUAL}.  Moreover, for any cylinder set $U$ we can find $(\alpha_k\in\ZZ)_{k\in K}$ and $(Y(k)\subset \ZZ^d)_{k\in K}$ 
 with $|K| < \infty$ such that
$$
\mu (U)=\sum_{k\in K}\alpha_k\left(\sum_{l=1}^\infty \mathbb{P}_{\xi_0=Y(k)}\{\xi_l=\emptyset\vert \xi_{l-1}\neq\emptyset \}\right).
$$
\end{thm} 
\begin{rem}
 Note that in some cases it is possible to exchange the role of the two states $\left(0 \leftrightarrow 1\right)$ in order 
to show ergodicity using the previous results.
\end{rem}
\begin{cor}\label{p=0}
Under the conditions in Theorem \ref{EXPR} we have that if  $\lambda(\emptyset)=0$ then the unique invariant measure is $\delta_0$, where $\delta_0$ is the Dirac measure on $^\infty 0^\infty$. Analogously, we have that if $\lambda(\emptyset)=1$ then the unique invariant measure  $\mu$ is $\delta_1$, where $\delta_1$ is the Dirac measure on $^\infty 1^\infty$.
\end{cor}
\begin{thm}\label{decay} 
Let $\eta_{.}$ be a one-dimensional probabilistic cellular automaton $\in \mathcal{C}$ of radius $r$ with $p(I_r)=:{\bf D}\in [0,1)$. Then, the unique 
invariant measure $\mu$ is shift-mixing. Also, if ${\bf D}\neq 0$, for any pair of cylinders $[U]_0=[u_0\ldots u_k]_0$, $[V]_0=[v_0\ldots v_{k'}]_0$ and $t\ge\vert U\vert+\vert V\vert$ we have 
\[  
\vert \mu ([U]_0\cap \sigma^{-t}[V]_0)-\mu ([U]_0)\times \mu ([V]_0)\vert \le \exp{(-a\times t)}\times K(U,V),
\]
where $\sigma$ is the shift on $\{0,1\}^{\ZZ}$, $a=1/2r\times\ln{(1/ {\bf D})}$ and $K(U,V)$ is a constant depending only on $U$, $V$, ${\bf D}$ and $r$. 
\end{thm}
\begin{rem}
This last result can be  extended to $d$-dimensional PCA.
\end{rem}
\subsection{Examples and comparison with known results}
\subsubsection{The Domany-Kinzel model}
This is a one-dimensional PCA $\eta_.$ of radius $r=1$ introduced in \cite{DK} with transition probabilities 
\begin{eqnarray*}
\mathbb{P}\{\eta_{t+1}(z)=1\vert\eta_{t}(z-1,z,z+1)=000 \mbox{ or } 010 \}&&=p(\emptyset )=p(\{0\})=a_0, \\
\mathbb{P}\{\eta_{t+1}(z)=1\vert\eta_{t}(z-1,z,z+1)=100 \mbox{ or } 110 \}&&=p(\{-1\}) \\
&&=p(\{-1,0\})=a_1, \\
\mathbb{P}\{\eta_{t+1}(z)=1\vert\eta_{t}(z-1,z,z+1)=001 \mbox{ or } 011 \}&&=p(\{1\})=p(\{0,1\})=a_1
\end{eqnarray*}
and
\begin{eqnarray*}
\mathbb{P}\{\eta_{t+1}(z)=1\vert\eta_{t}(z-1,z,z+1)=101 \mbox{ or } 111 \}&&=p(\{-1,1\}) \\
&&=p(\{-1,0,1\}) = a_2,
\end{eqnarray*}
where, for any subset $V \subset \mathbb{Z}, \eta(V) \in \{0,1\}^V$ denote the restriction of a configuration $\eta \in \{0,1\}^{\mathbb{Z}}$ to the set of positions in $V$.\\
Using Proposition \ref{cond}  we obtain that $\eta_.\in \mathcal{C}$ when  
$p(\{-1,1\})\ge p(\{-1\})+p(\{1\})-p(\emptyset)$, which is equivalent to the condition $a_2\ge 2a_1-a_0$.
From Theorem \ref{EXPR} the PCA $\eta_.$ is ergodic if $p(I_r)=p(\{-1,0,1\})=a_2<1$.
From Theorem \ref{decay} the unique invariant measure is shift-mixing 
with exponential decay of spatial correlation such that for any pair of cylinders $[U]_0$ and $[V]_0$ and for all 
$t\ge \vert U\vert +\vert V\vert$ we obtain 
$$
\vert \mu ([U]_0\cap \sigma^{-t}[V]_0)-\mu ([U]_0)\times\mu ([V]_0)\vert\le K\exp{(-(1/2\ln{(1/a_2)})t)},
$$
where $K$ can be explicitly computed (see the end of the Proof of Theorem \ref{decay}).
Using Theorem \ref{EXPR} we can compute, for example, the measure of the cylinder $[01]_0$ which 
is 
\begin{eqnarray*}
\mu ([01]_0)&=&\mu ([1]_1)-\mu ([11]_0)=\hat{\mu}(\{1\})-\hat{\mu}(\{0,1\}) \\
&=&\sum_{t=1}^{\infty}\mathbb{P}_{\xi_0=\{1\}}\{\xi_t=\emptyset\vert \xi_{t-1}\neq\emptyset\} \\
&+&\sum_{t=1}^{\infty}\mathbb{P}_{\xi_0=\{0,1\}}\{\xi_t=\emptyset\vert \xi_{t-1}\neq\emptyset \},
\end{eqnarray*}
where $\xi_.$ is the associated dual process.
\subsubsection{Two-dimensional example}
Let $\eta$ be a two-state, two-dimensional PCA of radius one.
In this case $I_1=\{(l,k)\vert -1\le l,k\le 1\}$ is a square of 9 sites.
The transition probabilities $\{p(J)\vert J\subset I_1\}$ of $\eta_.$ are defined by 
$p(J)=\alpha \sum_{k=0}^{\vert J\vert}C_k^9=\alpha\times 2^{\vert J\vert}$ 
where $C_k^l$ are the binomial coefficients.
 This PCA belongs to $\mathcal{C}$ 
 since for any $J\subset I_1$ we can write $\lambda (J)=\alpha$ and obtain that  $P(J)=\sum_{J'\subset J}\lambda (J')$.
This PCA is a kind of generalization to dimension 2 of  the Domany-Kinzel model 
(each site has the same weight) with only one parameter.
The sufficient ergodicity condition is $p(I_r)<1$ which implies that  
$\alpha\times 2^9<1$ $(\alpha<2^{-9})$ and the constant of decay of spatial correlation 
is $a=\frac{1}{2}\ln (1/(2^9\times\alpha))$.

\subsubsection{ Comparison with Dobrushin condition}

In \cite{DO}, Dobrushin gives sufficient ergodicity conditions for interacting particule systems.
Using our notation, these conditions applied to PCA can be translated as $\gamma <1$ (see \cite{MS91} and \cite{MS93}), where
$$\gamma=\sum_{j\in I_r}\sup_{J\subset I_r}\vert p(J\cup\{j\})-p(J) \vert.
$$ 
In the case of the Domany-Kinsel model, which belongs to the class $\mathcal{C}$, we obtain 
$\gamma=\sup_{J\subset I_r}\vert p(J\cup\{-1\})-p(J) \vert+ \sup_{J\subset I_r}\vert p(J\cup\{1\})-p(J) \vert
=2(a_2-a_1)$ since $\eta_.\in\mathcal{C}$ $(a_2\ge 2a_1-a_0)$.
If $a_2<1$ (condition of Theorem \ref{EXPR}) and $2(a_2-a_1)\ge 1$ the Dobrushin suficient conditions can not be applied.\\
For the two-dimensional example we have $\gamma =\alpha \left(\sum_{k=1}^{9}k\times C_k^9\right)$.
In this case $\gamma> p(I_r)$ and even if $\gamma<1$ the constant of decay of correlation 
$\frac{1}{2}\ln (1/(p(I_r))$ is greater than $\frac{1}{2}\ln (1/(\gamma ))$, 
the constant of decay of correlation given in \cite{MS91}.\\  
More generally, if a PCA belongs to $\mathcal{C}$ the sufficient condition $p(I_r)<1$ 
can be rewritten as  $p(I_r)=\sum_{J\subset I_r}\lambda(J)<1$ 
and the Dobrushin sufficient condition can be rewritten as 
$\gamma=\sum_{J\neq\emptyset,\,  J\subset I_r}\lambda(J)\times\vert J\vert<1$.
  

\section{Proof of Theorem \ref{EXPR} and Proposition \ref{cond}}
\subsection{PCA in $\mathcal{C}$ and their Dual Process}
 In \cite{LSS2006}, the authors give sufficient ergodicity conditions for one-dimensional multi-state PCA of radius one using a dual process satisfying equation \ref{ExtendedDUAL}. Here we will use an analogous dual process to give sufficient ergodicity conditions for two-state, $d$-dimensional PCA of radius $r$ using the following duality equation:

\begin{equation}\label{dual-eq}
\mathbb{E}_{\eta_0=x}[H(\eta_1,Y)]=\mathcal{D}(Y)\mathbb{E}_{\xi_0=Y}H[(x,\xi_1)],
\end{equation}
where $\eta_.$ is a PCA with state space
$\{0,1\}^{\ZZ^d}$. The state space of the dual process $\xi_{.}$ is the class of all finite subsets of $\ZZ^d$. 
As in \cite{LSS2006} we define the function $H$ by
 $$
 H(x,Y)=\left\{
 \begin{array}{l}
 1 \mbox{ , if } x(z)=1, \forall z\in Y\\
 0 \mbox{ , otherwise.}
 \end{array}
 \right.
 $$ 
The rule for the evolution of the process $\xi_t$ is given  by 
$$\xi_{t+1}=\displaystyle \cup_{z\in \xi_t} B(z)
$$
where  for any nonempty set $J\subset I_r$ we have 
$$
\mathbb{P}\Big[B(z)=\{z+j\vert j\in J\}\Big]=\pi(J)
$$ 
and 
$$
\mathbb{P}\big[B(z)=\emptyset\big]=\pi (\emptyset).
$$

Then, take the function $\mathcal{D}$ such that  $\mathcal{D}(Y)={\bf D}^{\vert Y\vert}$ for any finite subset $Y\subset\ZZ^d$, where ${\bf D}\in [0,1]$.
 Note that $\mathcal{D}(\emptyset )=1$ and $\emptyset$ is the unique absorbing state for this dual process.
\subsection{The functions $H$ and $\hat{\mu}$}
Note that, for this particular choice of $H$, we have
$$
\hat{\mu}(\ZZ^d)=\int_{X} H(x,\ZZ^d )d\mu (x) = \mu(^\infty 1^\infty) = 0
$$
and 
$$
\hat{\mu}(\emptyset) = \int_{X} H(x,\emptyset )d\mu (x) = \mu(\{0,1\}^{\ZZ^d}) = 1,
$$
where $X=\{0,1\}^{\ZZ^d}$ and $^\infty 1^\infty$ is the all one configuration $\left(1_{\mathbb{Z}^d}(x)\right)_{x \in \mathbb{Z}^d}$. 
The following Lemma is used in the proof of Theorem \ref{EXPR} 
and Theorem \ref{decay}. 

\begin{lem}\label{pro-dense}
The set of linear combinations of $\{H(.,y)\vert y\in\ZZ^d\}$ is dense in $C\left(\{0,1\}^{\ZZ^d},\RR\right)$, the set 
of continuous 
function from $\{0,1\}^{\ZZ^d}$ to $\RR$. For any cylinder $U=N(\Lambda,\varphi)\subset \{0,1\}^{\ZZ^d}$  (with $\Lambda\subset\ZZ^d, |\Lambda| < \infty$ and $\varphi\in A^{\ZZ^d}$) we have
$$
\mu (U)=\sum_{Y(i)} \alpha_i\hat{\mu}(Y(i)),
$$
where $\alpha_i\in\ZZ$, $Y(i)\subset\ZZ^d$ and $\max\{\vert Y(i)\vert\} < \infty$.
\end{lem}

\noindent {\bf Proof.} For the sake of simplicity, we only give the proof for the two-state, one-dimensional case.
The key point of the proof consists in showing that any cylinder $[U]_t:=[u_0\ldots u_n]_t$, 
($u_i\in \{0,1\}$ and $t,n\in\NN$) can be decomposed into a non-commutative sequence of subtractions and unions of 
intersections 
of cylinders of the type $[1]_t, t\in\ZZ$. We denote by $T([U]_t)$ this decomposition.
One way to accomplish this decomposition is to follow the following rules:
$$
 T([1]_t)=[1]_t, \hskip .5 cm  T([0]_t)=\{0,1\}^{\mathbb{Z}^d} - [1]_t.
$$
Then, for all $t,n\in\ZZ$ and $U=u_0\ldots u_n$ we have
$$
T([U1]_t)=T([U]_t)\cap [1]_{t+n+2}.
$$
Thus,
$$
T([U0]_t)=T([U]_t)-T([U]_t)\cap [1]_{t+2+n}.
$$
For instance,
$$
\begin{array}{ll}
T([100]_0)&=T([10]_0)-T([101]_0)\\
&=\left(T([1]_0)-T([11]_0)\right)-\left(T([10]_0)\cap [1]_2\right)\\
&=[1]_0-[11]_0-\left(([1]_0-[11]_0)\cap [1]_2\right)\\
&=\left([1]_0-[11]_0-\left([1]_0\cap [1]_2\right)\right)\cup [111]_0.\\
\end{array}
$$
Then, note that ${\bf 1}_{[1000]_0}$, the characteristic function of the cylinder $[1000]_0$,
can be written as 
\begin{eqnarray}
{\bf 1}_{[1000]_0}(x)&=&{\bf 1}_{[1]_0}(x)+{\bf 1}_{[111]_0}(x)-{\bf 1}_{[1]_0\cap [1]_2}(x)
-{\bf 1}_{[11]]_0}(x) \nonumber \\
&=& H(x,\{0\})+H(x,\{0,1,2\})-H(x,\{0,2\})-H(x,\{0,1\}). \nonumber
\end{eqnarray}
Since for any finite subset $Y\subset\ZZ$  we have ${\bf 1}_{\cap_{i\in Y} [1]_i}(x)= H(x,Y)$, it follows
that for all $n\in\NN$, $t\in\ZZ$ and $U\in {\{0,1\}}^n$ ${\bf 1}_{[U]_t}=\sum \alpha_i H(x,Y(i))$. This, in turn,
implies that the set of linear combinations of the set $\{H(.,Y)\vert Y\in \ZZ^d\}$ is dense in $C(\{0,1\}^{\ZZ^d})$.
We finish the proof by observing that for any cylinder $[U]_t$, we have 
\begin{eqnarray}
\mu ([U]_t)&=&\int {\bf 1}_{[U]_t}(x)d\mu (x) \nonumber \\
&=&\int \sum \alpha_i H(x,Y(i)) d \mu (x) \nonumber \\
&=& \sum\alpha_i\hat{\mu}(Y(i)). \nonumber
\end{eqnarray}

\hfill$\Box$

\begin{rem}
Using the definition of  $H$ taken in \cite{LSS2006} which takes into consideration the multi-state case, it is possible to prove Proposition \ref{pro-dense} for more general $d$-dimensional PCA.
\end{rem}
\subsection{ Proof of Theorem \ref{EXPR}}

We first establish a sequence of equalities between the transition probabilities 
of the PCA $(P(J)\vert J\in I_r)$ and the transition probabilities of the dual process 
($(\pi(J)\vert J\in I_r)$).

We can rewrite the right hand of equation (1.4) to  obtain  
$$
\mathbb{E}_{\eta_0=x}[H(\eta_1,Y)]=\mathbb{P}_{\eta_0=x}\{\eta_1(z)=1\; \forall z\in Y\}.
$$ 
Hence, using the independence property of $\eta_{.}$ we get that 
$$
\mathbb{P}_{\eta_0=x}\{\eta_1(z)=1\; \forall z\in Y\}=\displaystyle \prod_{z\in Y}\mathbb{P}_{\eta_0=x}\{\eta_1(z)=1\}.
$$ 
For the left hand of equation \ref{dual-eq}
we have 
$$
\mathbb{E}_{\xi_0=Y}[H(x,\xi_1)]=\mathbb{P}_{\xi_0=Y}\{x(z)=1\; , \forall z\in \xi_1\}.
$$
For any $x \in \{0,1\}^{\mathbb{Z}^d}$ we denote by $C_x$ the set $\{z\in\ZZ^d\vert x(z)=1\}$. Then
$$
\mathbb{P}_{\xi_0=Y}\{x(z)=1\; \forall z\in \xi_1\}=\mathbb{P}_{\xi_0=Y}\{\xi_1\subset C_x\}.
$$ 
Using the independence property of the dual process we can assert that 
$$
\mathbb{P}_{\xi_0=Y}\{\xi_1\subset C_x\}=\displaystyle \prod_{z\in Y}\mathbb{P}\{B(z)\subset C_x\}.
$$ 
Finally we can rewrite equation \ref{dual-eq} as  
\begin{equation}
\begin{array}{ll}
\displaystyle \prod_{z\in Y}\mathbb{P}_{\eta_0=x}\{\eta_1(z)=1\}&={\bf D}^{\vert Y\vert}
\displaystyle \prod_{z\in Y}\mathbb{P}\{B(z)\subset C_x\}\\
&=\displaystyle \prod_{z\in Y}{\bf D}\times\mathbb{P}\{B(z)\subset C_x\}
\end{array}
\end{equation}
which implies that 
\begin{equation}\label{eqgen}
\begin{array}{l}
\displaystyle \prod_{z\in Y}p(J_z)
=\displaystyle \prod_{z\in Y}{\bf D}\times \Delta_z,
\end{array}
\end{equation} 
where $J_z=\{i-z\vert i\in\{C_x\cap \{j+z\}\vert j\in I_r\}\}$ and $\Delta_z$ is given by
\begin{eqnarray}
\Delta_z &=& \pi(\emptyset ) + \sum_{i\in I_r}{\bf 1}_{1}(x(z+i))\times\pi(\{i\}) \nonumber \\
&+& \sum_{i,j\in I_r}{\bf 1}_{\{1\}}(x(z+i))\times{\bf 1}_{\{1\}}((x(z+j))\times\pi(\{i,j\})+\ldots \nonumber \\
&+& \sum_{i_1,\ldots , i_k\in I_r}\left(\prod_{l=1}^k{\bf 1}_{\{1\}}(x(z+i_k))\right)\times\pi(\{i_1,\ldots ,i_k\}) \nonumber \\
&+& \ldots + \left(\prod_{i\in I_r}{\bf 1}_{\{1\}}(x(z+i))\right)\times\pi(I_r). \nonumber
\end{eqnarray}

By simplicity of notation we  write $\pi(i_1,\ldots i_k)$ and $p(i_1,\ldots i_k)$ 
instead of $\pi(\{i_1,\ldots i_k\})$ and $p(\{i_1,\ldots i_k\})$.

Since equation \ref{eqgen} is true for all $x\in \{0,1\}^{\ZZ^d}$ 
 we obtain the following equations for $\pi (.)$,
 
\small
$$\begin{array}{ll}
p(\emptyset )&={\bf D}\pi(\emptyset )\\
p(i)&={\bf D}[\pi (\emptyset )+\pi(i)]\\
p(i,j)&={\bf D}[\pi (\emptyset )+\pi(i)+\pi(j)+\pi(i,j)] \\
p(i,j,k)&={\bf D}[\pi (\emptyset )+\pi(i)+\pi(j)+\pi(i,j)+\pi(i,k)+\pi(j,k)+\pi(i,j,k)]
\end{array}$$
\normalsize
where $i, j, k \in I_r$.\\
More generally, for any $0\le k\le \vert I_r\vert -1,$
\small
\begin{equation} \label{equationp}
p(i_0,... ,i_k)=
{\bf D}\Big[\pi (\emptyset )+\sum_{l=0}^k\pi(l)+\ldots +
\sum_{i=0}^{k-1}\sum_{l_0,... ,l_{i}\in \{i_0,\ldots ,i_k\}}\pi(l_0,... ,l_{i})
+\pi(l_0,l_1,... ,l_{k})\Big].
\end{equation}
\normalsize
\noindent Since
$$
\pi (\emptyset )+\sum_{k=0}^{\vert I_r\vert}\left(\sum_{l_0,l_1,\ldots l_k\in I_r}\pi(l_0,l_1,\ldots l_k)\right)
+\pi(I_r)=1,
$$
we get that ${\bf D}=p(I_r)$.
 
\bigskip
 

By definition, the dual process  is completely  
determined by the parameters $0\le \pi(J) \le 1$ ($J\subset I_r$).
From the sequence of equations \ref{equationp} the dual process  associated 
with the particular functions $H$ and $\mathcal{D}$ exists if the transition probabilities of the 
PCA satisfy $p(J)={\bf{D}}\sum_{J\subset I_r}\pi (J)$ with $0<p(I_r)\le I_r$.
In this case we have that $\lambda(J)={\bf{D}} \pi (J)$ and we claim that 
a PCA $\eta_.$ admits a dual process that satisfies the duality equation \ref{ExtendedDUAL} with  
particular functions $H$ and $\mathcal{D}$ given in section 4.1 
if and only if this PCA belongs to the class $\mathcal{C}$.

\medskip

\noindent
To show that the PCA is ergodic we need to verify the three conditions of Theorem \ref{MTHM}.\\
Condition $i)$ is verified since 
from Lemma \ref{pro-dense}, the set of linear combinations of functions belonging to 
$\{H(.,Y),\vert Y\in\ZZ^d\}$ is dense in $C\left(\{0,1\}^{\ZZ^d},\RR\right)$.\\ 
Condition $ii)$ is satisfied since $\sup_{Y\neq\emptyset} \{\mathcal{D}(Y)\}={\bf D}=p(I_r)<1$.\\ 
Condition $iii)$ follows from the fact that  
$H(.,\emptyset)=1$ and $\mathcal{D}(\emptyset)={\bf D}^{|\emptyset|}=1$.\\

\noindent
Since $\emptyset$ is the only absorbing state for $\xi_.$, 
 using  Theorem \ref{MTHM} (equation 1.2) 
we get that for any nonemptyset $Y\subset \ZZ^d$
$$
\hat{\mu}(Y)=\hat{\mu}(\emptyset) \mathbb{P}_{\xi_0=y}\{\xi_{\tau}=\emptyset\}
=\sum_{t=1}^{+\infty}\mathbb{P}_{\xi_0=y}\{\xi_{t}=\emptyset\vert \xi_{t-1}\neq\emptyset\}.
$$
From Lemma \ref{pro-dense}, for any cylinder set $U$ there exist $\alpha_k\in \RR$ and $Y(k)$ finite subset of 
 $\ZZ^d$ such that $\mu (U)=\sum \alpha_k\hat{\mu}(Y(k))$, which implies the last statement 
of Theorem \ref{EXPR}.

\hfill$\Box$
\bigskip

\subsubsection{Proof of Corollary  \ref{p=0}}

 When $\lambda (\emptyset)=1$, starting from any initial measure $\mu_0$, we obtain that $\mu_1=\delta_1$.
When $\lambda (\emptyset)=0$, Theorem \ref{EXPR} and Lemma \ref{pro-dense} imply that for each cylinder  
$U$ that does not contain the point $^\infty 0^\infty$ we get  
$$
\mu(U)=\sum\alpha_i\left(\sum_{k=0}^{\infty}\mathbb{P}_{Y_0=Y(i)}
\{Y_k=\emptyset\vert Y_{k-1}\neq\emptyset\}\right)=0
$$ 
since  $\pi(\emptyset)=\frac{\lambda (\emptyset)}{p(I_r)}=0$.
Finally we get that $\mu (^{\infty}0^{\infty})=1-\mu (\{0,1\}^{\ZZ^d}-^{\infty}0^{\infty})=1$ which 
finishes the proof.
\hfill$\Box$

\subsection{Proof of Proposition \ref{cond}}

Since the $\{\pi (J)\vert  J\subset I_r\}$ represent the transition probabilities of the 
dual process for all $J\in I_r$ one has $\pi(J)\ge 0$ and 
 Proposition \ref{cond} is a simple consequence of 
the following  Lemma .

\begin{lem}\label{lempi}
The transition probabilities $\pi()$ of the dual process satisfy
\small
 $$
\begin{array}{ll}
\pi (\emptyset )&=\frac{p(\emptyset)}{{\bf D}}\\
\pi(i)&=\frac{p(i)-p(\emptyset)}{{\bf D}}\\
\pi(i,j)&=\frac{1}{{\bf D}}[p(i,j)+p(\emptyset)-p(i)-p(j)]\\
\pi(i,j,k)&=\frac{1}{{\bf D}}[p(i,j,k)-p(\emptyset)+p(i)+p(j)+p(k)-p(i,j)-p(i,k)-p(j,k)]\\
\pi(i,j,k,l)&=\frac{1}{{\bf D}}[p(i,j,k,l)+p(\emptyset)-\sum_{l_0\in \{i,j,k,l\}} p(l_0)+
\sum_{\{l_0,l_1\}\subset\{i,j,k,l\}} p(l_0,l_1)\\
&-\sum_{\{l_0,l_1,l_2\}\subset \{i,j,k,l\}} p(l_0,l_1,l_2)]\\
\end{array}
$$
\normalsize
More generally, for any $0 \leq k \leq |I_r|-1$ and for any $j_0,\ldots ,j_k\in I_r$
\small
$$
\pi(j_0,\ldots ,j_k)=\frac{1}{{\bf D}}\left[(-1)^{k+1}p(\emptyset)+
\sum_{j=0}^{k} (-1)^{k-j}\sum_{\{l_0,\ldots ,l_j\}\subset \{j_0,\ldots ,j_k\}} p(l_0,\ldots ,l_j)\right].
$$
\normalsize
\end{lem}
{\bf Proof. of Lemma \ref{lempi}}
From the proof of Theorem \ref{EXPR} a PCA belongs to class $\mathcal{C}$ if and only if 
 the transitions probabilities $p()$ and $\pi()$ satisfy  the sequence of equations \ref{equationp}. 
We use  mathematical induction to solve the sequence of equations \ref{equationp}. 
For the two first iterations it is easily seen that   $\pi(\emptyset)=\frac{p(\emptyset)}{{\bf D}}$, $\pi(i)=
\frac{p(i)-p(0)}{{\bf D}}$ and $\pi(i,j)=\frac{1}{{\bf D}}[p(i,j)+p(0)-p(i)-p(j)]$.
 Then suppose that the order $k$ is true: 
$$
\pi(j_0,\ldots j_k)=\frac{1}{{\bf D}}\left[(-1)^{k+1}p(\emptyset )+
\sum_{j=0}^{k} (-1)^{k-j}\sum_{(l_0,\ldots ,l_j)\in \{j_0,\ldots ,j_k\}} p(l_0,\ldots l_j)\right].
$$ 
Using equation \ref{equationp} we obtain that  $\pi(j_0,\ldots ,j_{k+1})$ equals 
\begin{equation}\label{intermed}
\frac{1}{{\bf D}}\left[p(j_0,\ldots ,j_{k+1})-d\pi(\emptyset)- {\bf D} 
\sum_{j=0}^{k} \left(\sum_{(l_0,\ldots ,l_j)\in \{j_0,\ldots ,j_{k+1}\}} \pi(l_0,\ldots ,l_j)\right)  \right].
\end{equation}
Then we suppose the rank $k$ true and use  equation \ref{intermed} to obtain that 
 the term in $p(\emptyset )$ in  $\pi(j_0,\ldots ,j_{k+1})$  
is 
\begin{eqnarray}
&-&p(\emptyset )-\sum_{i=0}^{k}\left(\sum_{l_0,\ldots ,l_i\in\{j_0,\ldots ,j_{k+1}\}}(-1)^{i+1}p(\emptyset)\right) \nonumber \\
&=&p(\emptyset)\left(-1-\sum_{i=0}^k C_{i+1}^{k+2}(-1)^{i+1}\right) \nonumber \\
&=&p(\emptyset )\Big( -1+C_0^{k+2}(-1)^0+C_{k+2}^{k+2}(-1)^{k+2}-(1-1)^{k+2}\Big) \nonumber \\
&=&(-1)^{k+2}p(\emptyset ), \nonumber
\end{eqnarray}
where the constants $C_i^k$ represent the binomial coefficients.
Next we obtain that  the term in  $\sum_{l_0\in\{j_0,\ldots ,j_{k+1}\}}p(l_0)$ 
in $\pi(j_0,\ldots ,j_{k+1})$ is equal to 
\begin{eqnarray}
&-&\sum_{i=0}^k\sum_{(l_0,\ldots ,l_{i})\in\{j_0,\ldots ,j_{k+1}\}}
\left(\sum_{h_0\in\{l_0,\ldots ,l_{i}\}}p(h_0)\right)
(-1)^{i} \nonumber \\
&=&-\sum_{l_0\in\{j_0,\ldots ,j_{k+1}\}}p(l_0)\left(\sum_{i=0}^{k}C_{i}^{k+1}(-1)^{i}\right) \nonumber \\
&=&-\sum_{l_0\in\{j_0,\ldots ,j_{k+1}\}}p(l_0)\left((1-1)^{k+1}-C_{k+1}^{k+1}(-1)^{k+1}\right) \nonumber \\
&=& \sum_{l_0\in\{j_0,\ldots ,j_{k+1}\}}p(l_0)(-1)^{k+1}. \nonumber
\end{eqnarray}
Note that  $C_{i}^{k+1}$  represents  the number of ways to choose $l_1,\ldots ,l_{i}$ in $j_1,\ldots ,j_{k+1}$
when we have chosen $l_0$ and $j_0$.
More generally, for $0\le M\le k$,  the term in 
$\displaystyle \sum_{(l_0,\ldots ,l_M)\in\{j_0,\ldots ,j_{k+1}\}}$ $p(l_0,\ldots ,l_M)$ in $\pi(j_0,\ldots ,j_{k+1})$ is 
equal to
\begin{eqnarray}
&-&\sum_{i=M}^{k}\sum_{(l_0,\ldots ,l_{i})\in \{j_0,\ldots ,j_{k+1}\}}
\left(\sum_{(h_0,\ldots ,h_M)\in\{l_0,\ldots ,l_{j}\}}p(h_0,\ldots ,h_M)\right)(-1)^{i-M} \nonumber \\
&=&-\sum_{(l_0,\ldots ,l_M)\in\{j_0,\ldots ,j_{k+1}\}}p(l_0,\ldots ,l_M)
\left(\sum_{i=0}^{k-M} C_{i}^{k+1-M}(-1)^{i}\right) \nonumber \\
&=&-\sum_{(l_0,\ldots ,l_M)\in\{j_0,\ldots ,j_{k+1}\}}p(l_0,\ldots ,l_M)\left((1-1)^{k+1-M}-(-1)^{k+1-M}\right) \nonumber \\
&=&\sum_{(l_0,\ldots ,l_M)\in\{j_0,\ldots ,j_{k+1}\}}p(l_0,\ldots ,l_M)(-1)^{k+1-M}. \nonumber 
\end{eqnarray}

\hfill$\Box$

\section{Decay of Correlation}

For the sake of simplicity we study the decay of correlation for PCA with state space $\{0,1\}^{\ZZ}$.
An extension of this result to the multi-dimensional case is straightforward but requires too much notation. 

\subsection{Proof of Theorem \ref{decay}}

The proof of Theorem \ref{decay} requires the following two results. 
The second one is new and is a key point for the proof of Theorem \ref{decay}. 
The first one  seems to be well known. 
However, its proof can not be found or at least it is quite hard to be found so
 we provide a proof of that result.

Recall that $\mu$ stands for the unique invariant measure of an ergodic PCA.
\begin{pro}\label{si}
Every  invariant measure of an ergodic PCA  is shift-invariant.
\end{pro}

\begin{lem}\label{couplage}
Let 
$[U]_0$ and $[V]_0$ be two cylinders. If $\mu([U]_0)=\sum\alpha_i\hat{\mu}(A_i), 
\mu([V]_0)$ $=\sum\beta_i\hat{\mu}(B_i)$ and $t\ge \vert U\vert +\vert V\vert$, then
$$
\mu([U]_0\cap \sigma^{-t}[V]_0)=\mu([U]_0\cap [V]_t) = \sum\alpha_i\hat{\mu}(A_i)(*,t)\sum\beta_i\hat{\mu}(B_i),
$$
where 
$$\sum\alpha_i\hat{\mu}(A_i)(*,t)\sum\beta_i\hat{\mu}(B_i) := \sum_{i,j}\alpha_i\beta_j\hat{\mu}
(A_i\cup\{B_i+t\}).$$
\end{lem}

\subsection*{Proof of Theorem \ref{decay}}
If ${\bf D}=0$, then $p(\emptyset )=0$. From Corollary \ref{p=0}, $\mu=\delta_0$ and $\mu$ has exponential decay of 
correlation. For the remainder of this proof we therefore take $0<{\bf D}=p(I_r)<1$. 

For any finite subset $E$ of $\ZZ$ and $s\in\ZZ$, define $E+s :=
\{x+s : x\in E \}$.
We claim that for any finite subsets $E$ and $F$, if $t\ge 2Nr+\vert E\vert+\vert F\vert$ we have 
$$
\left\vert \hat{\mu}(E\cup \{F+t\})-\hat{\mu}(E)\times\hat{\mu}(F) 
\right\vert\le {\bf D}^{N+1}\frac{1}{1-{\bf D}}. 
$$
The proof of this claim uses Theorem 1 and  \ref{EXPR} which together say that 
 for any finite subset $E\subset\ZZ$, 
$\hat{\mu}(E)=\mathbb{P}_{\eta_0=E}\{\eta_{\tau}=\emptyset \}$. This, in turn, implies that
$$
\hat{\mu}(E)=\sum_{k=0}^{\infty}\mathbb{P}_{\eta_0=E}\{\tau =k\},
$$
where $\tau$ is the hitting time for the process $\eta_{.}$. In fact, by Lemma \ref{tau}, for any integer $N>0$ we have
$$
\left\vert \hat{\mu}(E)-\sum_{k=0}^{N}\mathbb{P}_{\eta_0=E}\{\tau =k\}\right\vert\le {\bf D}^{N+1}\frac{1}{1-{\bf D}}.
$$
Note that if  $s\ge 2ri+\vert E\vert +\vert F\vert $, where $i$ is any positive integer, then 
\begin{eqnarray}
\mathbb{P}_{\eta_0=E\cup \{F+s\}}\{\tau =i\}&=&\mathbb{P}_{\eta_0=E}\{\tau = i\}
\times\sum_{j=0}^i \mathbb{P}_{\eta_0=\{F+s\}}\{\tau = j\} \nonumber \\ 
&+& \mathbb{P}_{\eta_0=\{F+s\}}\{\tau =i\}\times\sum_{j=0}^i \mathbb{P}_{\eta_0=E\cup \{F+s\}}
\{\tau = j\}. \nonumber
\end{eqnarray}
It follows that  if $s\ge \vert E\vert +\vert F\vert +2N\times r$, then

$$
\sum_{i=0}^{N}\mathbb{P}_{E_0=E\cup \{F+s\}}\{\tau =i\}=
\sum_{i=0}^{N}\mathbb{P}_{E_0=E}\{\tau =i\}\times\sum_{i=0}^{N}\mathbb{P}_{E_0=\{F+s\}}\{\tau =i\}.
$$
This easily implies
\begin{equation}\label{muchap}
\left\vert \hat{\mu}(E\cup \{F+s\})-\hat{\mu}(E)\times\hat{\mu}(F) 
\right\vert\le {\bf D}^{N+1}\frac{1}{1-{\bf D}},
\end{equation}
for $s\ge \vert E\vert +\vert F\vert +2N\times r$, which proves our claim.

By Lemma \ref{pro-dense}, for any pair of cylinders 
$[U]_0$ and $[V]_0$, there exist finite sequences of sets $(A_i)$ and $(B_i)$ and finite sequences of real numbers
$\alpha_i$ and $\beta_i$ 
such that 
$$
\mu([U]_0)=\sum\alpha_i\hat{\mu}(A_i)$$
and 
$$
\mu([V]_0)=\sum\beta_i\hat{\mu}(B_i).
$$
Thus, by inequality \ref{muchap},
$$
\Big\vert \alpha_i\beta_j\hat{\mu}(A_i\cup\{B_i+s\}-\alpha_i\hat{\mu}(A_i)\times\beta_j\hat{\mu}(B_j))
\Big\vert \le \vert \alpha_i\beta_j\vert {\bf D}^{N+1}\frac{1}{1-{\bf D}}
$$
for any pair of subsets $A_i$ and $B_j$ of $\ZZ$ and for any $s\ge \vert  U\vert +\vert V\vert +2 N r$.

It follows from this that 
$$
\Big\vert \sum_{i,j}\alpha_i\beta_j\hat{\mu}\big(A_i\cup\{B_i+s\}\big)-
\sum_i\alpha_i\hat{\mu}(A_i)\times\sum_j\beta_j\hat{\mu}(B_j)\Big\vert 
\le F(U,V) {\bf D}^{N},
$$
where $F(U,V)=\sum_{i,j}\vert \alpha_i\beta_j\vert \frac{{\bf D}}{1-{\bf D}}$. 

Using Lemma \ref{couplage}, if 
$t\ge \vert U\vert +\vert V\vert$ 
we obtain 
$$
\Big\vert \mu\left([U]_0\cap \sigma^{-t}[V]_0\right)-\mu([U]_0)\times\mu([V]_0)\Big\vert \le K(U,V) 
\exp{\Big(-t\times\frac{\ln{(1/{\bf D})}}{2r}\Big)},
$$
where $K(U,V)=F(U,V) {\bf D}^{-(\frac{\vert U\vert+\vert V\vert}{2r})}$.  

Finally, it follows from Proposition \ref{si} that the invariant measure is shift-invariant and that the exponential decay of 
correlations of cylinders implies the mixing property.

\hfill$\Box$


\subsubsection{Proof of Proposition \ref{si}.}
It is sufficient to show that for any cylinder  $[U]_t$, where $U\in \{0,1\}^l$ for some $l\in\NN$ , we have 
$$
\mu (\sigma^{-1}[U]_t)=\mu ([U]_t).
$$
Since $\mu$ is the invariant measure of an ergodic PCA $\eta_{.}$,  there exits a sequence $(\mu_i)_{i\in\NN}$ which converges in the weak* topology to 
$\mu$, where $\mu_i$ is the distribution of a PCA $\eta_{.}$ at time $i$ starting from an initial distribution 
$\mu_0$.
It follows that for any cylinder $[U]_t$ we have 
$$\lim_{n\to\infty}\mu_n ([U]_t)=\mu ([U]_t).
$$
Since for any positive integer  $i$ we have 
$$
\mu_i([U]_t)=\sum_{V_j\in \{0,1\}^{n+1+2ir}} 
\mu_0 ([V_j]_{t-ir})\mathbb{P}_{\eta_0\in [V_j]_{t-ir}}\{\eta_i\in[U]_t\}, 
$$
 we can choose  $\mu_0$ as a shift-invariant probability measure. Hence, for any positive integer  $i$ 
and any cylinder $[U]_t$ we have   
$$
\mu_i ([U]_t)=\mu_i (\sigma^{-1}[U]_t),
$$
 which finishes the proof.

\hfill$\Box$

\subsubsection{Proof of Lemma \ref{couplage}} 
We prove the lemma using the principle of mathematical induction.
First we prove it for the case $\vert U\vert=1$ and $\vert V\vert =1$. Note that for any finite set 
$\{j_0,\ldots, j_k\}$, with $j_0,\ldots ,j_k\in\ZZ$,
\begin{eqnarray}  
\hat{\mu}(\{j_0,\ldots, j_k\})&=&\int_{\{0,1\}^\ZZ} H(\{j_0,\ldots, j_k\},x)d\mu(x) \nonumber \\
&=&\mu (\{\cap [1]_j\vert j\in\{j_0,\ldots, j_k\}\}) \nonumber
\end{eqnarray}
 and observe that for any $k\ge 1\in\NN$ we have 
 $$
 \mu([1]_0\cap [1]_k)=\hat{\mu}(\{0\}\cup\{k\}).
 $$
 Since $\mu ([0]_k)=1-\hat{\mu}(\{k\})=\hat{\mu}(\emptyset)-\hat{\mu}(\{k\})$ and  $\mu ([1]_0\cap [0]_k)
=\mu([1]_0)-\mu ([1]_0\cap [1]_k)=\hat{\mu}(\{0\})-\hat{\mu}(\{0,k\})$ we get, again, that
 $$
 \mu ([1]_0\cap [0]_k)=\hat{\mu}(\{0\}\cup\emptyset )-\hat{\mu}(\{0\}\cup\{k\}).
 $$
Furthermore, we have  $\mu ([0]_0\cap [1]_k)=\hat{\mu}(\{1\})-\hat{\mu}(\{0,k\})$. 

Finally, note that
\begin{eqnarray*}
\mu ([0]_0\cap [0]_k)&=&1-\mu ([1]_0\cap [1]_k)-\mu ([1]_0\cap [0]_k)-\mu ([0]_0\cap [1]_k) \\
&=& \hat{\mu}(\emptyset )-\hat{\mu}(\{0\}\cup \{k\})-\hat{\mu}(\{0\})+\hat{\mu}(\{0\}\cup \{k\})-\hat{\mu}(\{k\}) \\
&+& \hat{\mu}(\{0\}\cup \{k\}) \\
&=& \hat{\mu}(\emptyset )+\hat{\mu}(\{0\}\cup \{k\})-\hat{\mu}(\{0\})-\hat{\mu}(\{k\}) \\
&=& \left[\hat{\mu}(\emptyset )-\hat{\mu}(\{0\})\right](*,t)\left[\hat{\mu}(\emptyset )-\hat{\mu}(\{k\})\right],
\end{eqnarray*}
which finishes the proof in the case $\vert U\vert =\vert V\vert =1$.

\medskip

Now, suppose that $\mu([U]_0\cap \sigma^{-t}[V]_0)=\sum_{i,j}\alpha_i\beta_j\hat{\mu}(A_i\cup\{B_i+t\})$ 
is true for $\vert U\vert =\vert V\vert =n$. Consider $t\ge 2n+1$ and let $[U]_0, [V]_0$ be two cylinders such that 
$\mu ([U]_0)=\sum\alpha_i\hat{\mu}(A_i)$ and $\mu ([V]_0)=
\sum\beta_j\hat{\mu}(B_j)$. Since $\mu ([U1]_0)=\sum\alpha_i\hat{\mu}(A_i\cup \{\vert U\vert\})$, we get
\begin{eqnarray*}
\mu ([U1]_0[V1]_{t})&=&\sum_{i,j}\alpha_i\beta_j\hat{\mu}(A_i\cup \{B_j+t\}\cup\{\vert U\vert,\vert V\vert+t\}).
\end{eqnarray*} 
Noting that $\mu ([V1]_t)=\sum\alpha_i\hat{\mu}(\{B_i+t\}\cup \{\vert V\vert+t\})$ we get the desired result for 
the case $[U1]_0\cap [V1]_t$.

The result for the case $[U1]_0\cap [V0]_t$ follows by noting that
\begin{eqnarray*}
\mu ([U1]_0\cap [V0]_t)&=&\mu ([U1]_0\cap [V]_t)-\mu ([U1]_0\cap [V1]_t) \\
&=&\sum_{i,j}\alpha_i\beta_j\hat{\mu}(A_i\cup \{\vert U\vert\}\cup \{B_j+t\}) \\
&-&\sum_{i,j}\alpha_i\beta_j\hat{\mu}(A_i\cup \{B_j+t\}\cup\{\vert U\vert,\vert V\vert+t\}) \\
&=&\sum\alpha_i\hat{\mu}(A_i\cup \{\vert U\vert\})(*,t)\sum\alpha_i\hat{\mu}(\{B_i+t\}) \\
&-&\sum\alpha_i\hat{\mu}(\{B_i+t\}\cup \{\vert V\vert+t\}).
\end{eqnarray*}
It can also be shown that 
\begin{eqnarray}
\mu ([U0]_0\cap [V1]_t)&=&[\sum\alpha_i\hat{\mu}(A_i)
-\sum\alpha_i\hat{\mu}(A_i\cup \{\vert U\vert\})] \nonumber \\
&(*,t)&
[\sum\alpha_i\hat{\mu}(\{B_i+t\}\cup \{\vert V\vert+t\})]. \nonumber
\end{eqnarray}
Finally, using that 
\begin{eqnarray*}
\mu ([U0]_0\cap [V0]_t) &=& \mu ([U]_0\cap [V]_t)-\mu ([U1]_0\cap [V0]_t) \\
&-& \mu ([U0]_0\cap [V1]_t)- \mu ([U1]_0\cap [V1]_t)
\end{eqnarray*}
we can show  that 
\begin{eqnarray*}
\mu ([U0]_0\cap [V0]_t)&=&[\sum_{i,j}\alpha_i\beta_j\hat{\mu}(A_i\cup\{B_i+t\}) \\
&-&\sum_{i,j}\alpha_i\beta_j\hat{\mu}(A_i\cup\{B_i+t\}\cup\{\vert U\vert\})] (*,t) \\
&[&\sum_{i,j}\alpha_i\beta_j\hat{\mu}(A_i\cup\{B_i+t\}) \\
&-&\sum_{i,j}\alpha_i\beta_j\hat{\mu}(A_i\cup\{B_i+t\})\cup\{\vert V\vert +t\}]   
\end{eqnarray*}
which finishes the proof.

\hfill$\Box$



\section*{Final Questions}

\begin{itemize}
\item [i)] Is there exist an ergodic PCA such that the unique invariant measure is not shift-mixing?

\item [ii)] Is there exist an ergodic PCA  such that the  
invariant measure has non-exponential decay of correlation?



\end{itemize}

\begin{acknowledgements}
We would like to thank Marcelo Sobottka for a carefully reading of a previous version of this work and for many comments and suggestions. We also thank the referees for constructive remarks. During the realization of this work the first author received partial financial support from FAPESP,
grant 2006/54511-2. The second author received financial support from UFABC, grant 2008.
\end{acknowledgements}



\end{document}